\documentstyle[twoside,amssymb,12pt]{article}
\newcommand{\p}{\partial}

\newtheorem{prop}{Proposition}[section]

\newtheorem{theo}[prop]{Theorem}

\newtheorem{coro}[prop]{Corollary}
\newtheorem{lem}[prop]{Lemma}

\title{\sc On the case of Goryachev-Chaplygin
 and new examples of  integrable conservative
systems on $S^2$ } 

\author{{\sc Elena N. Selivanova}\thanks{Supported by DAAD.}\\ \\
\small  { Department of Geometry, Nizhny Novgorod State 
Pedagogical University}\\
\small  { 603000 Russia, Nizhny Novgorod, ul. Ulyanova 1}\\
\\  
\small  { Mathematisches Institut, Universit\"at T\"ubingen}   \\
\small  {Auf der Morgenstelle 10,  72076  T\"ubingen, Germany}   \\
\small  {e-mail: \ \ lena@moebius.mathematik.uni-tuebingen.de}\\ 
\small  { Fax: \ \ (07071)\ 294322}
}

\begin{document}
\date{}

\maketitle

\noindent
{\bf Abstract.} The aim of this paper is to describe a class of  conservative systems on
$S^2$ possessing an integral cubic in momenta.
We prove that this class of  systems consists off the case of Goryachev-Chaplygin, the one-parameter family 
of systems which has been found by the author in the previous paper and a  new 
 two-parameter family of conservative systems on $S^2$ possessing an integral cubic in momenta.

\thispagestyle{empty}

\section{\bf Introduction}

\noindent
Let $M$ be an $n$-dimensional  Riemannian manifold, and $U: M\to \bf R$ be a 
smooth function on $M$. For the Lagrangian $L: 
TM \to \bf R$
we consider $$L(\eta)=\frac{|\eta|^2}{2} + (U\circ\pi )(\eta)$$
where $\pi :TM\to M$ is the canonical projection, see \cite{Arn}. In local coordinates $q_1, \ ...\ ,  q_n$, $\dot q_1, \
...\ ,  \dot q_n$  on $TM$ we have
$$L(\eta)=\frac{1}{2}\sum {g_{ij}\dot q^i\dot q^j} + U(q).$$
Identifying $TM$ and $T^*M$ by means of the Riemannian metric, 
we get a Hamiltonian system
with the Hamiltonian $H: T^*M \to \bf R$ which 
 in local coordinates $q_1, \ ...\ , q_n$, $p_1, \ ... \ , p_n$ on $T^*M$ has the form
$$H=\frac{1}{2}\sum{g^{ij}p_ip_j} + U(q)=K+U.$$
We will call these Hamiltonian systems  {\it conservative systems} 
on $M$.

A smooth function $F : T^*M \to \bf R$ which is an integral of the 
Hamiltonian system with the Hamiltonian 
$H$ and which is independent of $H$ we will call an {\it integral} of this system 
of {\it degree} $m$ in momenta if in local coordinates $F$ has the form
$$F=\sum_{k_1+ \ ...\  + k_n\le m}{a_{k_1 \ ...\    k_n}(q)p_1^{k_1}  \ ...\   p_n^{k_n}}.$$ 

We will say that two Hamiltonians $H_1=K_1+U_1$ and $H_2=K_2+U_2$ 
are {\it equivalent} if there exists a diffeomorphism  $\phi$ of $M$ and a 
diffeomorphism
$\Phi$ of  $T^*M$
such that the diagram
$$\Phi: T^*M \to T^*M$$
$$\pi' \downarrow \ \ \ \ \ \ \ \ \ \ \downarrow \pi' $$
$$\phi: M\to M$$
is commutative, where   $\Phi$ is linear for   $p\in M$ fixed, and if
there are  some nonzero constants $\kappa_1$, $ \kappa_2$ such that
$\Phi^*(K_1)=\kappa_1 K_2$, $\phi^*(U_1)= \kappa_2 U_2$.

Clearly,  if the Hamiltonians $H_1$, $H_2$ are equivalent
and one of the corresponding systems possesses an integral of degree $m$ in momenta,
then the other system has the same property.
\medbreak
In this paper we consider integrable conservative systems on $S^2$ with the Hamiltonians of the form
\begin{equation}
H=\lambda(r^2)(r^2d\varphi^2 + dr^2) + \rho(r^2)r\cos\varphi=K +U
\label{formH}
\end{equation}
in polar coordinates $r, \varphi$ where $\lambda$, $\rho :{\bf R}^+\to {\bf R}$ are some smooth functions. 

\medbreak
There  is a known example of an integrable 
conservative system on $S^2$,  the case of Goryachev-Chaplygin 
in the dynamics of a rigid body, possessing an integral
 cubic in momenta. 
The total energy (the Hamiltonian) of this 
system has the following form
\begin{equation}
H=\frac{du_1^2 + du_2^2 + 4du_3^2}{4u_1^2 + 4u_2^2 + u_3^2} - u_1,
\label{gor}
\end{equation}
where $S^2$ is given by $ u_1^2 + u_2^2 + u_3^2=1$, (see 
 \cite{BKF}).

\medbreak
The aim of this paper is to classify (up to  equivalence, see above) 
the Hamiltonians 
of the form (\ref{formH}) of conservative systems on $S^2$, possessing an integral 
\begin{equation}
F=p_{\varphi}^3 + \kappa p_{\varphi}K + A
\label{formF}
\end{equation}
where $A$ is a linear polynomial in momenta and $\kappa$ is an arbitrary constant. 
We prove that these systems include the case of Goryachev-Chaplygin, a one-parameter family 
from \cite{Se} and a {\it new } two-parameter family of conservative systems on $S^2$ 
possessing an integral cubic in momenta.
It turns out that all these systems are due to smooth solutions of the following differential equation
\begin{equation}
x'x'''=xx'' - 2x''^2 + x'^2 + x^2.
\label{x}
\end{equation}
For this reason we start with investigations of this equation. We show
that there is a positive constant $T$ such that any smooth solution of (\ref{x})  
which exists everywhere can be obtained from the solutions $x=x(t): \bf R \to \bf R$ of 
the initial value problem  
\begin{equation}
x'x'''=xx'' - 2x''^2 + x'^2 + x^2, \ \  x(0)=0, \ x'(0)=1, \ x''(0)={\tau},
\label{ix}
\end{equation}
where $\tau\in [0, T]$,
or from $x(t)=\exp t$, $x(t)=\cosh t$ by a scaling $x \to \alpha_1 x$, $\alpha_1 - const$ or a linear translation of time $t\to \pm t + t_0$,
$t_0 - const$.

In the next section we obtain a criterion for integrability for geodesic flows of Riemannian
metrics as a partial differential equation. 
In the last section we consider a class of solutions of this equation and 
reduce the problem of classification of conservative systems on $S^2$ to the classification 
(up to a scaling or  linear translations) of smooth solutions of (\ref{x}).
We will show that 
the known case of Goryachev-Chaplygin is due to the solution of (\ref{ix}) if $\tau =T$ 
(if $\tau =0$ in (\ref{ix}) we get
 the  metric of constant positive curvature). Then we prove  the main theorem of 
classification and find new integrable cases on $S^2$. 

\section{\bf Smooth solutions } 

\medbreak
\begin{theo}
There is a positive constant $T$ such that any smooth solution of (\ref{x})  
which exists everywhere can be obtained  from $x(t)=\exp t$, $x(t)=\cosh t$ or 
from the solutions of 
the initial value problem  (\ref{ix}), where $\tau\in [0, T]$, 
by a scaling $x \to \alpha_1 x$, $\alpha_1 - const$ or a linear translation of time $t\to \pm t + \alpha _0$,
$\alpha_0 - const$.

For any solution $x_{\tau}$ of (\ref{ix}) with $\tau\in [0, T]$ there are smooth 
functions $\xi_{\tau}$, $\zeta_{\tau}$, $\mu_{\tau}$, $\nu_{\tau}$ such that
$$x_{\tau}'(t)=(\exp(-t))\xi_{\tau}(\exp(2t))=(\exp t)\zeta_{\tau}(\exp(-2t)),$$
$$(x_{\tau}'' - x_{\tau})x_{\tau}'^2(t)=(\exp t)\mu_{\tau}(\exp(2t))=(\exp(-t))\nu_{\tau}(\exp(-2t))$$
where $x_{\tau}'(t)$ is positive everywhere if and only if $\tau\in [0, T)$ and 
 there is one value of the parameter $t=t_0$
such that $x_{T}'(t_0)=0$.

\label{ex-th}
\end{theo} 

\noindent 
{\em Proof.} The initial value problem (\ref{x})
has a unique solution $x(t)=\Theta_{\tau}(t)$
which is positive on $(0,\varepsilon)$ and negative
on $(-\varepsilon, 0)$ for a sufficiently small 
$\varepsilon$.
\medbreak
\noindent
Let us consider the case $t>0$. 

The differential equation from  (\ref{x})
can be replaced by the following differential
equation of the second order 
$$\ddot{q}=\frac{1}{q}\left ( 1+2q^2-3q^4+\dot{q} -
7q^2\dot{q}-2{\dot{q}}^2\right) $$
with  $q(t)=R'(t)$ where $R(t)=\log x(t)$.
We may rewrite this equation as a system of differential equations of the first order:
\begin{equation}
\dot q=p, \ \
\dot p=\frac{1}{q}\left ( 1+2q^2-3q^4+p-7q^2p-2p^2 \right).
\label{syst1}
\end{equation}
Since  system (\ref{syst1}) is symmetric with respect to $q\mapsto -q$, 
$t\mapsto -t$,
 it  suffices to consider the case $q>0$.

In order to obtain the phase portrait of (\ref{syst1}) we may 
consider the following smooth system
\begin{equation}
\dot q=qp, \ \ 
\dot p=1+2q^2-3q^4+p-7q^2p-2p^2. 
\label{sms}
\end{equation}
The solutions of (\ref{syst1}) are obtained from the solutions of
 (\ref{sms}) by a
reparametrization. 
The system  (\ref{sms}) has four singular points: two saddle points 
$p=1, q=0$, $p=-\frac{1}{2}, q=0$
 and two nodes $p=0, q=\pm 1$. 
\medbreak
So, if a solution of (\ref{x}) is smooth everywhere, then it corresponds to a smooth solution
of (\ref{syst1}) and, therefore, it corresponds to $p=0, q=\pm 1$ or an orbit $p=-q^2+1, |q|<1$
or maybe two  orbits of (\ref{syst1}) 
where $q\to -\infty$ as $t\to 0-$ and $q\to +\infty$ as $t\to 0+$. 
In the cases $p=0, q=\pm 1$ and $p=-q^2+1, |q|<1$
 a solution of (\ref{x}) is obtained from $x(t)=\exp t$ and from $x(t)=\cosh t$
correspondingly by a scaling or a linear translation of time.
\medbreak
The aim of our further investigations is to show
that the orbits of (\ref{syst1}), corresponding to
the solutions of (\ref{x}), when $\tau$ belongs to a certain interval,
converge to the singular point $q=1, p=0$.
\medbreak
For any solution $x(t)=\Theta_{\tau}(t)$, $t\in (0, 
\varepsilon)$ of 
(\ref{x}) there is an orbit $\{ \Gamma_{\tau}: p=p_{\tau}(q)\}$ 
in the phase spase of (\ref{syst1}).

For any solution of (\ref{x}) it holds 
$R(t)\to -\infty$ as 
$t\to 0+$ and, therefore, $q(t)=R'(t)\to + \infty$,
 $p(t)=q'(t)\to -\infty$ as $t\to 0+$.
Thus $p_{\tau}(q)\to -\infty$ as $q\to +\infty$.

We may now consider only the orbits of (\ref{syst1}) where 
$p\to -\infty$ as $q\to +\infty$.

Show that if $\tau_1>\tau_2$, then $p_{\tau_1}(q)>p_{\tau_2}(q)$.
Indeed, for a solution $x(t)=\Theta_{\tau_1}(t)$ of (\ref{x}) it holds
\begin{equation}
\tau_1=\lim_{t\to 0+}{ \frac{\Theta''_{\tau_1}(t)}{\Theta'_{\tau_1}(t)} }=
\lim_{q\to +\infty}{\frac{q^2 + p_{\tau_1}(q)}{q}}.
\label{ta_u}
\end{equation}
Note that the function $x(t)=\sinh t$ satisfies (\ref{x}) when $\tau=0$.
The related orbit of (\ref{syst1}) has then the following form 
$\{\Gamma_0: p=p_0(q)=-q^2 + 1\}$.

So, the orbits of (\ref{syst1}), corresponding to the solutions of (\ref{x}),
for $\tau \ge 0$ converge to the singular point $q=1, p=0$. 
\medbreak
We prove below that orbit $(*)$, where 
$p\to -\frac{1}{2}$ as $q\to 0+$,
 corresponds to the solution
of (\ref{x}) when $\tau$ is equal to a negative constant $-T$ and all
orbits of (\ref{syst1}), lying between $(*)$ and $\Gamma_0$ correspond
to the solutions of (\ref{x}) for $\tau\in (-T,0)$.
\medbreak
Assume first that there exists a constant $\tau_0<0$ 
such that  orbit $\Gamma_{\tau_0}$ does not converge
to the point $q=1, p=0$.

Consider the set $W$ of the orbits of (\ref{syst1}), 
lying between  $\Gamma_{\tau_0}$ and $\Gamma_0$.
Show that for any orbit in $W$ the value $q$ becomes infinite in a 
finite time interval.

In fact, for any solution of (\ref{syst1}) it holds  
\begin{equation}
\int^{q(t)}\frac{dq}{p}=t + const.
\label{oc}
\end{equation}
For orbits in $W$ we have $p<-q^2 +1$. Hence, the left hand side of (\ref{oc})
is bounded for any orbit in $W$ as $q\to +\infty$.
Without loss of generality we assume that for these 
 solutions of (\ref{syst1})  $t$ vanishes  as $q$ becomes infinite. 

We conclude that  for any orbit from $W$ it holds by (\ref{ta_u}) 
$$(\tau_0+o(1))q<p+q^2<1$$
and for a corresponding solution of (\ref{x}) we get for $t \to 0+$
$$\tau_0+o(1)<\frac{x''(t)}{x'(t)}<1.$$
Thus, for a corresponding solution of (\ref{x}) the function  $(\log x'(t))'$ is 
bounded in an interval, containing $0+$.
 Hence,  $\lim_{t\to 0+}{x'(t)}$
is bounded  for a solution $x(t)$ of (\ref{x}), 
corresponding to an orbit in $W$ and, therefore, $x(t)=\frac{x'(t)}{q}\to 0$
as  $t\to 0+$. 
 
Let us consider the orbit $(*)$  where $p=p^*(q)$.
If $(*)$ is the same as $\Gamma_{\tau_0}$ then $T=-\tau_0$. 

If $(*)$ is not the same as $\Gamma_{\tau_0}$ then it belongs to $W$ and,
hence, corresponds to the solution of (\ref{x})
for $\tau=-T$, where $T$ is equal to a positive constant  and, moreover, 
$T<-\tau_0$. 

Thus we have shown that for any $\tau\in (-T, T)$ where $T$ is a
 positive constant the system (\ref{x}) has a solution for all $t\ge 0$.

\medbreak
Let us assume now that there is no such orbit $\Gamma_{\tau_0}$ and, so, 
all orbits of (\ref{syst1}), corresponding to the solutions of (\ref{x}),
 converge to the singular point $q=1, p=0$.

Clearly, the orbit $(*)$ corresponds to a solution $x(t)=\eta(t)$ of 
the differential equation in (\ref{x}).
As  mentioned above, the function $q(t)=\frac{\eta'(t)}{\eta(t)}$ 
becomes infinite in a finite  time interval of $t$, say, as $t\to 0+$.
Since $p^*(q)<-q^2 - kq$ for any $k$ and
 for  sufficiently large $q$, for $\eta(t)$ it holds
\begin{equation}
\lim_{t\to 0+}{\frac{{\eta''(t)}}{{\eta'(t)}}}=-\infty.
\label{p}
\end{equation}
Note that $\eta(0)$ is finite because $R(0)-R(t)=\int_{t}^{0}{q(t)dt}<0$ for $t>0$ and 
$\eta(t)=\exp R(t)$. 
So, there are two cases: $\eta(0)=0$ or $\eta(0)\ne 0$.

Assume  that $\eta(0)\ne 0$. It follows that $\eta'(0+)=+\infty$ 
and from (\ref{p}) we get $\eta''(0+)=-\infty$.

Rewrite now the differential equation in (\ref{x})  in the following form 
$$\frac{x'''}{x}=-3\frac{x''-x}{x'} -2\frac{(x''-x)^2}{xx'} + \frac{x'}{x}.$$
Therefore,
it holds $$\frac{\eta'''}{\eta}=-3\frac{p^*(q) +q^2 -1}{q} -2\frac{(p^*(q) + q^2 -1)^2}{q} + q.$$
We obtain   $\eta'''(0+)=-\infty$. This is a contradiction  $\eta''(0+)=-\infty$.

Assume  that  $\eta(0)=0$ and  $\eta'(0)\ne 0$.
 From (\ref {p}) we get $\eta''(0+)=-\infty$. 
Rewriting  the differential equation in (\ref{x})  in the following form
$$x'''=x'-3\frac{(x''-x)x}{x'} - 2\frac{(x''-x)^2}{x'},$$
we obtain again  $\eta'''(0+)=-\infty$.

Therefore, we only have to consider one case  $\eta(0)=0$ and  $\eta'(0)=0$.
Taking into account that ${\eta'}(t)=\eta(t)q(t)$ and $q(t)>0$ as $t\to 0+$, we 
conclude that ${\eta'}(t)\to 0+$ as $t\to 0+$ but on the other hand  from 
(\ref {p})
it follows that ${\eta''(t)}<0$.

These contradictions finally show  that there is an orbit $\Gamma_{\tau_0}$
which does not converge to the point $q=1, p=0$
and, therefore, for any $\tau\in (-T, T)$, where  
 \begin{equation}
T=-\lim_{q\to +\infty}{\frac{p^*(q)+q^2}{q}}<\infty,
\label{T}
\end{equation}
solutions of (\ref{x}) exist on $[0, +\infty)$.

Since a solution $x(t)=\Theta_{\tau}(t)$ of (\ref{x}) equals  $-\Theta_{-\tau}(-t)$ if 
$t\le 0$, solutions of (\ref{x}) exist on $(-\infty, +\infty)$.
\medbreak
We show now the asymptotic behaviour of the solutions of 
(\ref{ix}) with $\tau\in [0, T]$ where $T$ is defined by (\ref{T}).

Put $s=\exp(-2y)$ and 
$$g(s)=\sqrt{s}x\left(-\frac{1}{2}\log s\right).$$
Then the initial value problem (\ref{x}) can be rewritten as 
 \begin{equation}
g'''=3\frac{g'' g'}{g-2g's} + s\frac{g''^2}{g-2g's},
 \  g(1)=0, g'(1)=-\frac{1}{2}, g''(1)=\tau 
\label{g}
\end{equation}
We will consider $\tau\in (-T, T]$.

Compute now for a solution of (\ref{x}):
\begin{equation}
x'(t)=(\exp t)(g-2g's),
\label{d}
\end{equation}
\begin{equation}
x''(t)-x(t)=4\exp {(-3t)}g'' (s).
\label{dd}
\end{equation}

\noindent
Let us consider the  system (\ref{syst1}).
Since the eigenvalues of the Jacobian of (\ref{syst1}) at $(q,p)=(1, 0)$ 
  are equal to  $-2$,  $-4$, there exist 
functions $P$ and $Q$ of class $C^1$ such that 
 $$q=Q(s)=1+Cs+o_1(s),$$ $$p=P(s)=-2Cs+o_2(s),$$
where $C$ is a constant, see \cite{Har}.

We may write:
$$R(t)=\int^{t}q=t+\psi_1(s).$$
Let us show that $\psi_1$ is of class  $C^1$. By differentiation we obtain
$$\frac{d\psi_1(s)}{ds}=-\frac{q-1}{2s}=-\frac{Cs+o_1(s)}{2s}\in C^{0}.$$
Write now for a solution of (\ref{x})
$$x(t)=\exp {R(t)}=(\exp t)\exp {\psi_1(s)}=(\exp t)g(s),$$
where we have used that $g=r^{-1}u=\exp(-t)u=\exp(-t)x(t)$,
and, therefore, $g(0)=\exp {\psi_1(0)}\ne 0$. 
\medbreak
\noindent
For a solution of (\ref{x}) which can be extended into infinity it
 holds  $$x'(t)=x(t)q(t)=(\exp t)(\exp {\psi_1(s)})(1+Cs+o_1(s))= $$ $$
=(\exp t) g(s)(1+Cs+o_1(s)), s=\exp (-2t).$$ 
On the other hand, from (\ref{d}) we obtain 
$$g(s)(1+Cs+o_1(s))=g(s)-2sg'(s). $$ 
Thus, for any solution of  (\ref{g}) where $\tau\in (-T, T]$
it holds $g'(0)<\infty$.
\medbreak
\noindent
As mentioned above, the solution of (\ref{x}) for $\tau=0$ has the form $x(t)=\sinh t$.

Let us show that for  any solution of  (\ref{g}) where $\tau\in (-T, T]$ 
and $\tau\ne 0$  it holds $0<|g''(0)|<\infty$.

We compute
$$x''(t)-x(t)=x(t)(p+q^2-1)=(\exp t)(\exp {\psi_1(s)})(P(s)+Q^2(s)-1)=$$
$$=(\exp t)(\exp {\psi_1(s)})(-2Cs+o_2(s)+(1+Cs+o_1(s))^2-1).$$
Thus,
$$x''(t)-x(t)=(\exp t)\psi_2(s),$$ where
$\psi_2\in C^1$, $\psi_2(s)=o_3(s)$.
Rewrite the differential equation in (\ref{x}) in the following form
$$(x''' - x')x'=-3(x'' - x)x - 2(x'' -x)^2.$$
We obtain  either $x''(t)-x(t)\equiv 0$ (then $\tau=0$)   
or
$$(\log |x''(t)-x(t)|)'=-3\frac{1}{q}-2|x''(t)-x(t)|(\exp t\exp {\psi_1(s)}(1+Cs+o_1(s)))^{-1}=$$
$$=-3\frac{1}{q}-2{|\psi_2(s)|}
(\exp {\psi_1(s)}(1+Cs+o_1(s)))^{-1}.$$ 
So, $$(\log |x''(t)-x(t)|)'= -3(1-Cs)+o_4(s).$$ 
Then it follows 
$$\left(\log \left((\exp t)|\psi_2(s)|\right)\right)'=-3(1-Cs)+o_4(s).$$ Thus,
$$(t+\log |\psi_2(s)|)'=-3(1-Cs)+o_4(s).$$
 By integrating we get
$$(\log |\psi_2(s)|)=-4t-\frac{3C}{2}s+\psi_3(s),$$ where
$\psi_3\in C^1$.

Thus,
$$ |\psi_2(s)|=(\exp(-4t))\exp ({-\frac{3C}{2}s+\psi_3(s)})=
s^2\psi_4(s),$$ since $s=\exp(-2t)$. So,
$\psi_4 \in C^1$ and  $\psi_4(0)=const \ne 0$. Taking into account
(\ref{dd}), we obtain
$|g''(0)|=\frac{1}{4}\psi_4(0)$ and therefore $0<|g''(0)|<\infty$.

So, we conclude that 
the solutions of (\ref{g}) for $\tau\in (-T, T]$ 
are of class  $C^{\infty}$ in zero.
\medbreak
\noindent
Now the functions $\zeta_{\tau}$, $\nu_{\tau}$ can be calculated in terms of $g$.
From (\ref{d}) $\zeta_{\tau}(0)>0$ in view of $g(0)>0$. Therefore $\zeta_{\tau}>0$ everywhere on $[0, +\infty)$.

Since a solution $x(t)=x_{\tau}(t)$ of (\ref{x}) equals  $-x_{-\tau}(-t)$ 
 for $\tau\in (-T, T)$ and 
   $x_{T}(t+t_0)=x_{T}(-t+t_0)$, the theorem follows. 

\hfill  $\Box$ 

\section{\bf A criterion for  integrability}

\noindent 
Consider a metric $ds^2=\Theta(u, v)(du^2 + dv^2)$ in conformal
coordinates $u, v$. It can  also be written as
\begin{equation}
ds^2=\theta(w, \bar w)dwd\bar w
\label{met}
\end{equation}
 where $w=u+ i v$. 
The geodesic flow of $ds^2$ is a Hamiltonian system with Hamiltonian
\begin{equation}
H=\frac{p_wp_{\bar w}}{{4\theta}(w, \bar w)}.
\label{ham}
\end{equation}
A polynomial $F$ in momenta $p_u$, $p_v$ 
can be also written as $$F=\sum_{k=0}^{n}{b_k(w,\bar w)p_w^kp_{\bar w}^{n-k}}$$
where $$b_k=\overline {b_{n-k}}, \ \ k=0,...,n.$$ 

If the polynomial $F$ is an additional integral of the geodesic flow with 
 the Hamiltonian (\ref{ham}), then $\{F,H\}=0$ and the following holds
\begin{equation}
\theta \frac{\partial b_{k-1}}{\partial w}+(n-(k-1))b_{k-1}
\frac{\partial \theta}{\partial w}+
\theta \frac{\partial b_{k}}{\partial \bar w}+kb_{k}
\frac{\partial \theta}{\partial \bar w}=0,
\label{systpde}
\end{equation}
where $k=0,... ,n+1 $ and $ b_{-1}=b_{n+1}=0$. 
Substituting   $k=0$ and $k=n+1$ in (\ref{systpde}) we get  immediately 
$$  \frac{\partial}{\partial \bar w}b_0\equiv 0$$ and 
$$\frac{\partial}{\partial  w}b_{n+1}\equiv 0.$$

One may show that if a polynomial in momenta integral $F$  
of the geodesic flow of 
(\ref{met})  is independent of the Hamiltonian 
and an integral of smaller degree, then there is a conformal coordinate system $z=z(w)$ 
of this  metric such that the coefficients of $F$ for $p_z$ and $p_{\bar z}$
 are equal  to $1$ identically.

\begin{theo} Let $ds^2=\lambda(z, \bar z)dzd\bar z$ ($z=\varphi + iy$) be a metric
such that there exists a function $f:{\bf R}^2\mapsto \bf R$, satisfying the 
following conditions 
$$
\lambda=\frac{1}{4}\left(\frac{\p ^2 f}{\p \varphi^2}+\frac{\p ^2 f}{\p y^2}\right),$$
and
\begin{equation}
\frac{\p}{\p \varphi}\left(\left(\frac{\p ^2 f}{\p \varphi^2}-\frac{\p ^2 f}{\p y^2}\right)\left( 
\frac{\p ^2 f}{\p \varphi^2}+\frac{\p ^2 f}{\p y^2}\right)\right)=2
\frac{\p}{\p y}\left(\left(\frac{\p ^2 f}{\p \varphi\p y}\right)\left( 
\frac{\p ^2 f}{\p \varphi^2}+\frac{\p ^2 f}{\p y^2}\right)\right)
\label{eqpde}
\end{equation}
Then the geodesic flow of $ds^2$ possesses an integral cubic in momenta.
\medbreak
If the geodesic flow of a metric $ds^2$ possesses an integral 
which is cubic in momenta and it does not depend on 
the Hamiltonian and an  integral of smaller degree then there exist
  conformal coordinates $\varphi, y$ and a function   $f:{\bf R}^2\mapsto \bf R$ such that 
$ds^2=\lambda(z, \bar z)dzd\bar z$ where $z=\varphi +iy$ and (\ref{eqpde}) holds.
\label{loc-th}
\end{theo} 

\noindent
{\em Proof.}  We will consider
integrals  $$F=\sum_{k=0}^{n}a_k(z, \bar z)p_z^kp_{\bar z}^{n-k}$$ 
of the geodesic flow of $ds^2$ where $a_0=a_n\equiv 1$. 

 The system (\ref{systpde})  has  then the following form
\begin{equation}
3\frac{\p \lambda}{\p z}+\frac{\p (a_1\lambda)}{\p \bar z}=0,
\label{1.7}
\end{equation}
\begin{equation}
\frac{\p (a_1\lambda^2)}{\p z}+\frac{\p (a_2\lambda)}{\p \bar z}=0.
\label{1.8}
\end{equation}

 (\ref{1.7}) holds if and only if there exists a function $h$ such that
$$\lambda=\frac{\p h}{\p \bar z}\ \mbox{and} \ \ a_1\lambda=-3\frac{\p h}{\p  z}.$$
(\ref{1.8}) can be also rewritten in the following form
$$\frac{\p (a_1\lambda^2)}{\p z}=-\frac{\p (\overline {a_1}\lambda^2)}{\p \bar z}=-
\overline{\frac{\p (a_1\lambda^2)}{\p z}}.$$
Therefore, the system (\ref{1.7}), (\ref{1.8}) is equivalent to the following condition
$$\mbox{Re}\ \frac{\p (a_1\lambda^2)}{\p z}=\mbox{Re}\ \frac{\p }{\p z}\left(\frac{\p h}{\p z}
 \frac{\p h}{\p \bar z}\right)=0$$
where $$\mbox{Im}\  \frac{\p h}{\p \bar z}=\mbox{Im}\  \lambda(z, \bar z)=0,$$
i.e. there is a real function $f$ such that $\frac{\p f}{\p  z}=h$ and  (\ref{eqpde}) holds.

\hfill  $\Box$ 

The equation (\ref{eqpde}) has been obtained first in \cite{Hall}.
We need to give hier another proof because further we need the following 
corollaries which are due to the above proof of Theorem \ref{loc-th}.  

\begin{coro}
If in a conform coordinate system $\varphi, y$ of $ds^2=\lambda(\varphi,y)(d\varphi^2  + dy^2)$ 
there is an integral of the geodesic flow
of $ds^2$  which
has the form
\begin{equation}
F=\alpha_3(p_{z}^3 + a_1p_{z}^2p_{\bar z} + 
\overline {a_1}p_{z}p_{\bar z}^2+ p_{\bar z}^3), \ z=\varphi+iy, \ \alpha_3 - const \ne 0,
\label{F'}
\end{equation}
then there is a function $f$ such that (\ref{eqpde}) holds.
\label{F}
\end{coro}

\begin{coro}
If in a conform coordinate system $\varphi, y$ of $ds^2=\lambda(\varphi,y)(d\varphi^2  + dy^2)$ there 
is a function $f$ such that (\ref{eqpde}) holds, then there is an integral $F$ which
has the form
$$
F=p_{z}^3 + a_1p_{z}^2p_{\bar z} + 
\overline {a_1}p_{z}p_{\bar z}^2+ p_{\bar z}^3, \ z=\varphi+iy.
$$
\label{c2F}
\end{coro}

\section {The Classification}

\begin{theo}
Let 
\begin{equation}
ds^2=(\Psi_1(r^2)r\cos\varphi + \Psi_2(r^2))(r^2d\varphi^2 + dr^2)
\label{metric}
\end{equation}
be a metric on $S^2$ in polar coordinates. Then the geodesic flow of $ds^2$
possesses an integral of the form
\begin{equation}
F=iw^3p_{w}^3 + b_1p_{w}^2p_{\bar w} + \overline {b_1}p_{w}p_{\bar w}^2 - 
i\bar w^3p_{\bar w}^3, \ \ w=r\cos\varphi + ir\sin\varphi
\label{metricF}
\end{equation}
if and only if 
$$\Psi_1(r^2)r=((\psi''-\psi)(\log r))r^{-2}$$ and
$$\Psi_2(r^2)=cr^{-2}\psi'^{-2}(\log r)$$ or
$$\Psi_2(r^2)=ar^{-2}\frac{\psi'^2 - \psi^2 + b}{\psi'^2}(\log r) $$ 
where $\psi$ is a smooth solution of 
(\ref{x}) which exist everywhere on $(-\infty, +\infty)$ and $a$, $c$, $p$ are appropriate constants.
\label{geo}
\end{theo} 

\noindent
{\em Proof.} Assume that the geodesic flow of $ds^2$
possesses an integral of the form (\ref{metricF}).
Let us consider $z=\varphi + iy=-i\ln w$. Then (\ref{metric}) can be rewritten
 in the following form
$$ds^2=(\Psi_1(\exp(2y))(\exp y)\cos \varphi + \Psi_2(\exp(2y)))(\exp(2y))
(d\varphi^2 + dy^2).$$ Thus, 
\begin{equation}
ds^2=(\Psi_3(y)\cos \varphi + \Psi_4(y))(d\varphi^2 + dy^2)=\lambda(\varphi, y)
(d\varphi^2 + dy^2).
\label{metric'}
\end{equation}
for some functions $\Psi_3$, $\Psi_4$.
An integral of the form (\ref{metricF}) has then the form (\ref{F'}).
Now we can apply Corollary \ref{F}. So, there is 
a function $f$ such that (\ref{eqpde}) holds.

Let us consider a function $h(\varphi, y)=\psi(y)\cos\varphi + \xi(y)$
where $\psi'' - \psi=\Psi_3$ and $\xi''=\Psi_4$. Then a function $f$ 
from (\ref{eqpde}) has the form $f=h+\delta$ where $\delta_{\varphi\varphi}
+\delta_{yy}=0$.
Then for the function $a_1$ from (\ref{F'}) we get
\begin{equation}
a_1=-3\frac{f_{\varphi\varphi}-f_{yy}}{f_{\varphi\varphi}+f_{yy}} +
6i\frac{f_{\varphi y}}{f_{\varphi\varphi}+f_{yy}}=
-3\frac{h_{\varphi\varphi}-h_{yy}+\delta_{\varphi\varphi}-
\delta_{yy}}{\lambda(\varphi, y)}+6i\frac{h_{\varphi y}+
\delta_{\varphi y}}{\lambda(\varphi, y)}.
\label{a_1}
\end{equation}
Thus, the functions $\delta_{\varphi y}$ and  $\delta_{\varphi \varphi}-
\delta_{yy}$ are periodic in $\varphi$. Therefore,
$$\delta(\varphi, y)=P_1(\varphi, y)
 + a(\varphi^2-y^2) + b\varphi y +\sum_{k\in \bf N}
{(C_k(\exp(-ky))+D_k(\exp(ky)))\sin k\varphi }$$
\begin{equation}
+ \sum_{k\in \bf N}{ (E_k(\exp(-ky))+G_k(\exp(ky)))\cos k\varphi} 
\label{delta}
\end{equation}
where $a$, $b$, $C_k$, $D_k$, $E_k$, $G_k$ are some constant and 
$P_1(\varphi, y)$ is a linear polynomial in $\varphi$, $y$.

For the function $b_1$ from (\ref{metricF}) we have 
$b_1(w, \bar w)=a_1|w|^2w$ is bounded if $|w|=0$, 
since $ds^2$ is a metric on $S^2$. 
Then from (\ref{a_1}) and (\ref{delta})
we obtain $C_k$, $D_k$, $E_k$, $G_k$ are equal to zero when $k\ge 2$.

Thus, we have to consider functions $f$ of the form 
\begin{equation}
f(\varphi, y)=\psi(y)\cos\varphi  + \xi(y) + a(\varphi^2 - y^2)
+(C_1(\exp(-y))+D_1(\exp y))\sin\varphi
\label{f()}
\end{equation}
satisfying (\ref{eqpde}). Taking into account that $ds^2$ is a metric on
$S^2$ from (\ref{eqpde}) for (\ref{f()}) we obtain $b=0$ and $C_1=D_1=0$.
 
So, we must consider solutions of  (\ref{eqpde}) of the following form
$$f(\varphi, y)=\psi(y)\cos\varphi  + \xi(y) +
 a(\varphi^2 - y^2), \ a - const.$$
Then (\ref{eqpde}) is equivalent to the following conditions 
$$\lambda=(\psi''-\psi)\cos x + \xi''$$ where 
\begin{equation}
\xi''\psi'' + (\xi''\psi')'=2a(\psi''-\psi) 
\label{8}
\end{equation}
and 
\begin{equation}
\psi'\psi'''=\psi\psi''-2\psi''^2+\psi'^2+\psi^2.
\label{9}
\end{equation}

Multiply (\ref{8}) by $\psi'$:
$$2\xi''\psi'\psi'' + \xi'''\psi'^2=2a(\psi''-\psi)\psi'$$
and integrate
$$\xi''(y)=a\frac{\psi'^2(y) - \psi^2(y) + b}{\psi'^2(y)}, \ b - const \ \  \mbox{if} \ \ a\ne 0$$
or $$\xi''(y)=c\psi'^{-2}(y),  \ c - const \ \ \mbox{if} \ \ a=0.$$
\bigbreak

Conversely integrability of the corresponding geodesic flows follows immediately from 
Corollary \ref{c2F} and above calculations.

\hfill$\Box$

Thus, due to the Maupertuis's principle we obtain the following corollary from Theorem \ref{geo}.

\begin{coro}
Hamiltonian systems  with the Hamiltonians of the form
\begin{equation}
H=\frac{d\varphi^2 + dy^2}{\psi'^2(y)} -
 (\psi''(y)-\psi(y))\psi'^2(y)\cos \varphi
\label{H}
\end{equation}
and 
\begin{equation}
H_b=\frac{\psi'^2(y)-\psi^2(y)+b}{\psi'^2(y)}\left(d\varphi^2 + dy^2\right) -
 \frac{\psi'^2(y)(\psi''(y)-\psi(y))}{\psi'^2(y)-\psi^2(y)+b}\cos \varphi,
\label{H_p}
\end{equation}
where   $\psi$ is a solution of (\ref{x}),
possessing an integral cubic in momenta of the form (\ref{formF}).
\label{cons}
\end{coro}

Now we prove the main theorem.

\begin{theo} A Hamiltonian system with 
the Hamiltonian 
of the form (\ref{formH}) is a conservative system on $S^2$, possessing an integral 
of the form (\ref{formF}) and  not possessing an integral quadratic or linear in momenta 
if and only if the corresponding Hamiltonian is equivalent to one of the following

1. (\ref{H}) where $\psi$ is a solution of (\ref{ix}) for $\tau\in (0, T)$.

2. (\ref{H_p}) where $\psi$ is a solution of (\ref{ix}) for $\tau\in (0, T)$ and 
$$ b>b_*(\tau)=\max_{y\in \bf R} {(\psi^2 - \psi'^2)} \ \ \mbox{or} \ \ 
b<b^*(\tau)=\min_{y\in \bf R} {(\psi^2 - \psi'^2)} $$
where $b_*(\tau)<+\infty$ and $b^*(\tau)>-\infty$ for $\tau\in (0, T)$. 

3. (\ref{H_p}) where $\psi$ is a solution of (\ref{ix}) for $\tau=T$ and $p=\psi^2(y_0)$,
when $\psi'(y_0)=0$, which is in fact the case of Goryachev-Chaplygin in the dynamics of a
rigid body.
\label{main-th}
\end{theo}

\noindent
{\em Proof.} From Theorem \ref{geo} and the Maupertuis's principle we know that 
if a Hamiltonian system with 
the Hamiltonian 
of the form (\ref{formH}) defines a conservative system on $S^2$, possessing an integral 
of the form (\ref{formF}), then (\ref{formH}) is equivalent to (\ref{H})
or (\ref{H_p}) where $\psi$ is a smooth solution of (\ref{x}) which  exists everywhere
on $(-\infty, +\infty)$. So, we must only describe when in this case 
(\ref{H}), (\ref{H_p}) define a conservative system on $S^2$.
From  Theorem \ref{ex-th} it follows that we must consider only 
solutions of (\ref{ix}) for $\tau\in [0, T]$ and $\psi(y)=\exp y$, $\psi(y)=\cosh y$.
We note that if $\psi$ is a solution of (\ref{ix}) for $\tau=0$ (then $\psi(y)=\sinh y$) 
or  $\psi(y)=\exp y$,  $\psi(y)=\cosh y$, then the Hamiltonians (\ref{H}), (\ref{H_p})
define simply a metric of constant curvature where there are two independent linear integrals.

So, we must only consider solutions $\psi$  of (\ref{ix}) when $\tau\in (0,T]$.

Let us write then  
(\ref{H}), (\ref{H_p}) in polar coordinates $\varphi$, $r=\log y$ and $\tilde\varphi=-\varphi$,
$\tilde r=-\log y$. Using  Theorem \ref{ex-th} we compute
$$H=\frac{1}{\xi_{\tau}^2(r^2)}(r^2d\varphi^2 + dr^2) - \mu_{\tau}(r^2)r\cos \varphi $$
$$=\frac{1}{\zeta_{\tau}^2(\tilde r^2)}
(\tilde r^2d\tilde \varphi^2 + d\tilde r^2) - 
\nu_{\tau}(\tilde r^2)\tilde r\cos \tilde \varphi $$
and
$$H_b=\frac{\Phi_{\tau}(r^2)+b+1}{\xi_{\tau}^2(r^2)}(r^2d\varphi^2 + dr^2) - 
\frac{\mu_{\tau}(r^2)r\cos \varphi}{\Phi_{\tau}(r^2)+b+1}$$
$$=\frac{-\tilde \Phi_{\tau}(\tilde r^2)+b+1}{\zeta_{\tau}^2(\tilde r^2)}
(\tilde r^2d\tilde \varphi^2 + d\tilde r^2) -\frac{ 
\nu_{\tau}(\tilde r^2)\tilde r\cos \tilde \varphi}{-\tilde \Phi_{\tau}(\tilde r^2)+b+1} $$
where $$\Phi_{\tau}(t)=\int_{1}^{t}{\mu_{\tau}(s)\xi_{\tau}^{-1}(s)ds},$$
$$\tilde \Phi_{\tau}(t)=\int_{1}^{t}{\nu_{\tau}(s)\zeta_{\tau}^{-1}(s)ds}.$$

Consider $\tau=T$. In Theorem \ref{ex-th} it has been proved that for the 
corresponding solution $\psi$  of (\ref{ix}) there is $y=y_0$ such that $\psi'(y_0)=0$, and
therefore there is $r=r_0>0$ such that $\xi_T(r_0^2)=0$. Thus, in the case $\tau=T$ 
(\ref{H}) does not define
 a system on $S^2$ and  (\ref{H_p}) in this case defines a system on $S^2$ if and only if
$p=\psi^2(y_0)$.

Consider $0<\tau<T$. Then (\ref{H}) defines a family of conservative systems on $S^2$ which 
has been found in \cite{Se}. In \cite{Se} it has been proved that for $\tau_1\ne \tau_2$
the corresponding Hamiltonians are not equivalent, the systems possess nontrivial cubic integrals
(there is no quadratic or linear integral) and no Hamiltonian from this family is equivalent to 
the case of Goryachev-Chaplygin.

Let us find the admissible values of the parameter $b$ in this case. From the above expressions for $H_b$ 
in polar coordinates we have 
$$b+1>\max_{[0,1]}{-\Phi_{\tau}}=M_1(\tau) \ \ \mbox{and} \ \ b+1>\max_{[0,1]}{\tilde \Phi_{\tau}}=\tilde M_1(\tau)$$
or
$$b+1<\min_{[0,1]}{-\Phi_{\tau}}=M_2(\tau) \ \ \mbox{and} \ \ b+1<\min_{[0,1]}{\tilde \Phi_{\tau}}=\tilde M_2(\tau).$$ 
So, $b>b_*(\tau)=\max {\{ M_1(\tau), \tilde M_1(\tau) \}}-1$ or 
$b<b^*(\tau)=\min {\{ M_2(\tau), \tilde M_2(\tau) \}}-1$ where $b_*(\tau)<+\infty$ and $b^*(\tau)>-\infty$ 
for $\tau\in (0, T)$.

With the same arguments as in \cite{Se} we can prove that the corresponding systems do not possess  integrals
linear or quadratic in momenta. 

Let us show that the Hamiltonians 
(\ref{H_p}) where $\psi$ is a solution of (\ref{ix}) for $\tau\in (0, T)$ 
are not equivalent for different values of
the parameters $\tau$ and $p$. 

Assume that the Hamiltoians of the form (\ref{H_p}) for $\tau_1$, $b_1$ and $\tau_2, b_2$ are equivalent.

We will use the following lemma which 
is in fact a particular case of 
some results of Kolokol'tsov, 
published in his Ph.D. Dissertation,
(Moscow State University, 1984), another proof one can find in \cite{Se}.

\begin{lem} Let 
\begin{equation}
ds^2=\lambda(r^2)(r^2d\varphi^2 + dr^2)
\label{polya}
\end{equation}
be a metric on $S^2$ in polar coordinates $r, \varphi$. 
Then $ds^2$ can be written in the form (\ref{polya}) in polar coordinates 
$\tilde r, \tilde \varphi$ if 
and only if $\tilde r=Dr^{\pm 1}$, $D - const$ and $\tilde \varphi=\pm \varphi +\varphi_0$,
$\varphi_0 - const$.
\label{L}
\end{lem}

Thus, from Lemma \ref{L} it follows that there are some constants $C_0\ne 0$, $C_3\ne 0$, $y_1$ such that
\begin{equation}
\psi_1''(y) - \psi_1(y)=C_0(\psi_2''(y+y_1) - \psi_2(y+y_1)),
\label{V}
\end{equation}
\begin{equation}
\frac{\psi_1'^2(y) - \psi_1^2(y) + b_1}{\psi_1'^2(y)}=C_3\frac{\psi_2'^2(y+y_1) - \psi_2^2(y+y_1) + b_2}
{\psi_2'^2(y+y_1)}
\label{K}
\end{equation}
where $\psi_1$, $\psi_2$ are solutions of (\ref{ix}) for $\tau=\tau_1$ and $\tau=\tau_2$ correspondingly.
So, we get from (\ref{V})
$$\psi_1(y)=C_0\psi_2(y+y_1) + C_1\exp y + C_2\exp(-y)$$
for some constants $C_1$, $C_2$. 
As in the proof of Theorem \ref{ex-th} let us write 
for a solution $\psi$ of (\ref{x}): $\psi(y)=(\exp y)g(\exp (-2y))$. It has been shown that
$g$ satisfies then the differential equation from (\ref{g}). 
So, from our assupmtion it follows that there are two solutions of  the differential equation from (\ref{g})
$g_1$ and $g_2$ such that 
$g_1(s)=g_2(s)+C_1 + C_2s$. Then by substituting  $g_1$  in the differential equation from (\ref{g})
we obtain $C_1=C_2=0$. Thus, $$\psi_1(y)=C_0\psi_2(y+y_1).$$ From the initial conditions of (\ref{ix})
we get $0=\psi_1(0)=C_0\psi_2(y_1)$ and therefore $\psi_2(y_1)=0$ but $\psi_2(0)=0$ and $\psi_2'(y)$ is positive
everywhere. So, $y_1=0$ and from $\psi_1'(0)=\psi_2'(0)=1$ we obtain $C_0=1$. Thus, $\tau_1=\tau_2$.

Prove now that $b_1=b_2$. So, we have from (\ref{K})
$$\psi_1'^2(y) - \psi_1^2(y) + b_1=C_3(\psi_1'^2(y) - \psi_1^2(y) + b_2),$$ 
and therefore, either $C_3=1$ and $b_1=b_2$ or $\psi_1'^2(y)-\psi_1^2(y)\equiv const$
 that is not true for $\tau\ne 0$. 
So we get $b_1=b_2$.
\medbreak

Assume that a Hamiltonian from   $(1.)$ is equivalent to a Hamiltonian from $(2.)$. Then we obtain 
in the same way as above that for a solution $\psi_{\tau}$ of  (\ref{ix}) it holds 
$\psi_{\tau}'^2(y) - \psi^2(y)\equiv const$ and therefore $\tau=0$. Thus,  no Hamiltonian from   $(1.)$
is equivalent to a Hamiltonian from $(2.)$.

\medbreak
Let us show that no Hamiltonian from this family 
is equivalent to the case of Gorychev-Chaplygin. There are polar coordinates $r, \varphi$ 
such that (\ref{gor}) can be rewritten as 
\begin{equation}
H=\gamma_1(r^2)(r^2d\varphi^2 + dr^2) - \gamma_2(r^2)r\cos\varphi=
\gamma_1(\tilde r^2)(\tilde r^2d\varphi^2 + d\tilde r^2) - 
\gamma_2(\tilde r^2)\tilde r\cos\varphi
\label{gor2}
\end{equation}
where $\tilde r=r^{-1}$. 

Assume that the Hamiltonian (\ref{H_p}) where $\psi$ is a solution of (\ref{ix}) for
 $\tau=\tau_0\in (0, T)$ and $b=b_1$ is equivalent
to the Hamiltonian of the case of Gorychev-Chaplygin. This means also from the symmetry of (\ref{gor2}) 
that (\ref{V}) and (\ref{K})
hold for $b_1=b_2$ and $\tau_1=\tau_0$, $\tau_2=-\tau_0$. In the same way as above we obtain $\tau_0=0$ but then
the Hamiltonian of the case of Gorychev-Chaplygin is equivalent to the Hamiltonian of a metric of constant 
positive curvature. 
\medbreak
So, we described all conservative systems on $S^2$ with the Hamiltonians 
of the form (\ref{formH}) possessing an integral 
of the form (\ref{formF}) and  not possessing an integral quadratic or linear in momenta.
We proved also that no  Hamiltonian from  $(1.)$ or $(2.)$ is equivalent to 
the Hamiltonian of the case of Goryachev-Chaplygin. On the other hand the case of Goryachev-Chaplygin belongs to
the  class of conservative systems on $S^2$ with the Hamiltonians 
of the form (\ref{formH}) possessing an integral 
of the form (\ref{formF}) and  not possessing an integral quadratic or linear in momenta, see \cite{Se}. 
So, the Hamiltonian of the case of Goryachev-Chaplygin is equivalent to 
the Hamiltonian (\ref{H_p}) where $\psi$ is a solution of (\ref{ix}) for $\tau=T$ and $b=\psi^2(y_0)$,
when $\psi'(y_0)=0$.

\hfill$\Box$

\bigbreak

\end{document}